\documentclass{amsart}
\usepackage{amssymb,amsmath,amscd, dsfont}
\usepackage[all]{xy}
\newcommand{\commentout}[1]{}
\newtheorem{lem}{Lemma}[section]
\newtheorem{thm}{Theorem}[section]
\newtheorem{cor}{Corollary}[section]

\newcommand{\mc}[1]{\mathcal#1}
\newcommand{\f}[1]{{\mathfrak{#1}}}

\newcommand{\mb}[1]{\mathbb#1}
\begin{document}
\pagenumbering{arabic}
\title{Unitary Representations and Heisenberg Parabolic Subgroup}
\author{Hongyu HE}\footnote{This research is partially supported by the NSF grant DMS 0700809 and by LSU.}
\footnote{keywords: Parabolic Subgroups, Heisenberg group, Stone-Von Neumann theorem, Mackey Analysis, branching Formula, unitary representations, Kirillov Conjecture, Symplectic Group}
\address{Department of Mathematics, Louisiana State University, Baton Rouge, LA 70803}
\email{hongyu@math.lsu.edu}
\date{}
\begin{abstract}
In this paper, we study the restriction of an irreducible unitary representation $\pi$ of the universal covering
$\widetilde{Sp}_{2n}(\mb R)$ to a Heisenberg maximal parabolic group $\tilde P$.  We prove that if $\pi|_{\tilde P}$ is irreducible, then $\pi$ must be a highest weight module or a lowest weight module.  This is in sharp constrast with the $GL_n(\mathbb R)$ case. In addition, we show that for a unitary highest or lowest weight module, $\pi|_{\tilde P}$ decomposes discretely. We also treat the groups $U(p,q)$ and $O^*(2n)$. 
 \end{abstract}
\maketitle
\section{Introduction}
Let $\mathbb F= \mathbb C, \mathbb R$. Let $GL_n(\mathbb F)$ be the general linear group on $\mathbb F^n$. Let $P_1$ be the maximal parabolic subgroup preserving a one dimensional subspace in $\mathbb F^n$. Let $\pi$ be an irreducible unitary representation of $GL_n(\mathbb F)$. Consider the restriction of $\pi$ onto $P_1$. Kirillov conjectured that $\pi|_{P_1}$ is irreducible. Kirillov's conjecture was proved by Sahi using Vogan's classification (~\cite{sah} ~\cite{vogan}). Recently, Baruch established
Kirillov's conjecture without Vogan's classification (~\cite{baru}). \\
\\
Generally speaking, for other semisimple Lie groups $G$, the restriction of an irreducible unitary representation of $G$ to a maximal parabolic subgroup is hardly irreducible. Nevertheless, as proved by Howe and Li, for irreducible low rank representations, their restrictions to a certain maximal parabolic subgroup remain irreducible (~\cite{howes} ~\cite{li}). In this situation, the restriction uniquely determines the original representation. However, it is not clear whether there are other representations whose restriction to a fixed maximal parabolic group is irreducible. \\
\\
Now Let $G=Sp_{2n}(\mb R)$ and {\bf $n \geq 2$}. Let $P$ be the maximal parabolic subgroup that preserves a one-dimensional isotropic subspace of the symplectic space $\mathbb R^{2n}$. Decompose the identity component of $P$  as $Sp_{2n-2}(\mb R) A H_{n-1}$ where $H_{n-1}$ is the Heisenberg group and $A \cong \mathbb R^{+}$. We call $P$ a Heisenberg parabolic subgroup of $G$. Let $\tilde G$ be the universal covering of $G$. Let $Mp_{2n}(\mb R)$ be the unique double covering of $G$. For any subgroup $H$ of $G$, let $\tilde H$ be the preimage of $H$ under the universal covering.  The classification of irreducible unitary representation of $\tilde P$ can be obtained directly by Mackey analysis. \\
\\
As an example, take the linear group $P$. Let $C(H_{n-1})$ be the center of $H_{n-1}$. Let $\pi$ be an irreducible unitary representation of $P$. If $\pi|_{C(H_{n-1})}$ is trivial, then $\pi$ is in one-to-one correspondence with irreducible unitary representations of  maximal parabolic subgroup of $Sp_{2n-2}(\mb R)$ with levi factor $GL_1(\mb R) Sp_{2n-4}(\mb R)$. Suppose $\pi|_{C(H_{n-1})}$ is not trivial. Let $\rho_{\pm}$ be the two irreducible unitary representation of dilated Heisenberg group $A H_{n-1}$. Then $\hat{P}$ is parametrized by a triple $(\rho_{\pm}, \tau, \pm)$ where $\tau$ is a genuine irreducible unitary representation of $Mp_{2n-2}(\mb R)$ and $\pm$ corresponds to the two representations of the component group of $P$. Extend $\rho_{\pm}$ to a unitary representation of $Mp_{2n-2}(\mb R) A H_{n-1}$, and extend $\tau$  trivially to $A H_{n-1}$. Every irreducible unitary representation of $P_0$ can thus be written as $\rho_{\pm} \otimes \tau $. Moreover, $\rho_{\pm} \otimes \tau $ can be extended to an irreducible unitary representation of $P$. So Every irreducible unitary representation of $P_0$ can be written as $\rho_{\pm} \otimes \tau 
\otimes \mathbb C_{\pm}$.\\
\\
In this paper, all tensor product of Hilbert spaces will mean the completion of the algebraic tensor product. All Hilbert spaces are assumed to be separable. We use $\pi$ to denote both the representation and the underlying Hilbert space.\\
\\
For simplicity, let us absorb the parameter $\pm$ into $\tau$. Any unitary representation $\pi$ of $P$ can then be written as
$$[\rho_{+} \otimes \tau_{+}] \oplus [\rho_{-} \otimes \tau_{-}] \oplus \tau_{0},$$
here $\tau_{0}|_{C(H_{n-1})}$ is trivial. Hence every irreducible unitary representation of $G$ can also be written in this form.  \\
\\
Notice that $\rho_{+}|_{Mp_{2n-2}(\mb R)}$ is equivalent to $\omega(n-1) \otimes \mathbb C^{\infty}$ where $\omega(n-1)$ is the oscillator representation of $Mp_{2n-2}(\mb R)$ and $\mathbb C^{\infty}$ is an infinite dimensional trivial representation of $Mp_{2n-2}(\mb R)$. $\rho_{-}|_{Mp_{2n-2}(\mb R)}$ is equivalent to  $\omega(n-1)^* \otimes \mathbb C^{\infty}$.

\begin{thm}[See ~\cite{howe}]
Let $\pi$ be a nontrivial  irreducible unitary representation of $\widetilde{Sp}_{2n}(\mb R)$. Let $\tilde P_0$ be the identity component of $\tilde P$. Then there are two unitary representations $\tau_{+}(\pi)$ and $\tau_{-}(\pi)$ of $\widetilde{Sp}_{2n-2}(\mb R)$ such that
$$\pi_{\tilde P_0} \cong [\rho_{+} \otimes \tau_{+}(\pi) ] \oplus [\rho_{-} \otimes \tau_{-}(\pi)],$$
$$\pi|_{\widetilde{Sp}_{2n-2}(\mb R)} \cong [ \omega(n-1) \otimes \tau_{+}(\pi)  \oplus \omega^*(n-1) \otimes \tau_{-}(\pi)] \otimes \mathbb C^{\infty}.$$
In the first identity, $\tau_{\pm}(\pi)$ extends trivially to representations of $\tilde P_0$.
\end{thm}
This theorem is established by Howe for the double covering $Mp_{2n}(\mb R)$ (\cite{howe}). Howe's argument essentially extends to the universal covering of $Sp_{2n}(\mb R)$. \\
\\
If $\pi$ is a unitarily induced representation from a unitary representation of $P$, then $\tau_{+}$ and $\tau_{-}$ are quite easy to compute. The issue of computing the map $\pi \rightarrow \tau_{\pm}(\pi)$ for smaller representations is rather complex. For the two constituents of the oscillator representation, $\omega(n)_{\pm}$, $\tau_{+}$ is trivial and  $\tau_{-}$ is zero. Moreover, if we take into consideration of the two connected component of $P$, $\tau_{+}$ for $\omega(n)_{+}$ will be the trivial representation and $\tau_{+}$ for $\omega(n)_{-}$ will be the sign character. \\
\\
By Mackey analysis, $\pi|_{\tilde P}$ is irreducible, if and only if one of $\tau_{\pm}(\pi)$ vanishes and the other is irreducible. Our main result is the following.

\begin{thm} Let $\pi$ be an irreducible unitary representation of $\tilde G$. If $\pi|_{\tilde P}$ is irreducible, then $\pi$ must be either a highest weight module or lowest weight module. In addition, for $\pi$ a unitary lowest weight module, $\pi|_{\tilde P} \cong \rho_{+} \otimes \tau_{+}(\pi)$ where $\tau_+(\pi)$ decompose discretely into a direct sum of lowest weight modules of $\widetilde{Sp}_{2n-2}(\mb R)$; for $\pi$ a unitary highest weight module, $\pi|_{\tilde P} \cong \rho_{-} \otimes \tau_{-}(\pi)$ where $\tau_-(\pi)$ decompose discretely into a direct sum of highest weight modules of $\widetilde{Sp}_{2n-2}(\mb R)$.
\end{thm}

It is not clear whether $\tau_{\pm}(\pi)$ is irreducible for $\pi$ a highest or lowest weight module. For some highest (lowest) weight modules, $\tau_{\pm}(\pi)$ is irreducible. In fact, for  $\widetilde{Sp}_2(\mb R)$, $\tau_{\pm}(\pi)$ will always be irreducible. For $n \geq 2$, decomposing $\pi|_{\tilde P}$ is quite difficult. Generally speaking 
$\pi|_{\tilde P}$ does not decompose according to the $K$-types. \\
\commentout
{
The map $\pi \rightarrow \tau_{\pm}(\pi)$ has some nice properties. For example, R. Howe shows that if $\pi$ is of pure rank $l$, then $\tau_{\pm}(\pi)$ is of pure rank $l-1$ (Cor. 2.12 \cite{howe}). Roughly, this puts a bound on the matrix coefficients of $\tau_{\pm}(\pi)$ from above. In the following theorem, we put a bound on the matrix coefficients of $\tau_{\pm}(\pi)$ from below. 
\begin{thm} Let $\pi$ be an irreducible unitary representation of $\widetilde{Sp}_{2n}(\mb R)$. Suppose that for any $\tilde U(n-1)$-finite vectors $u, v$ in $\tau_{\pm}(\pi)$, 
$$(\tau_{\pm}(\pi)(\exp H) u, v) \leq \| u \| \| v \| \exp \lambda(H)$$
where $H$ is in the positive Weyl chamber $\f a_{n-1}^+$ of the maximal split Cartan $\f a_{n-1}$. Then
$$(\lambda-\bold{\frac{1}{2}})(H) \geq  \mu(H)  \qquad (H \in \f a_{n-1}^+).$$
Here $\mu$ is a leading exponent of $\pi$ (\cite{knapp}).
\end{thm}
}
\\
In this paper, we derive some equivalent conditions for $\pi$ being a nontrivial highest weight module. One of the condition can be stated as follows.

\begin{thm} Let $Sp_2(\mb R)$ be a subgroup of $G$ that fixes a nondegenerate $2n-2$ dimensional symplectic subspace. Let $N$ be a unipotent subgroup in $Sp_2(\mb R)$. Identify $\hat N$ with the real line. Then $\pi$ is a nontrivial irreducible unitary highest or lowest weight module if and only if $\pi|_{N}$ is supported on half of the real line.
\end{thm}

In this paper, we also treat the groups $U(p,q)$ and $O^*(2n)$. The group $P$ will be a maximal subgroup whose nilradical is a Heisenberg group. We call such $P$ a Heisenberg parabolic subgroup. The detailed results are stated in Theorems ~\ref{main1u}, ~\ref{main2u}, ~\ref{main1o}, ~\ref{main2o}.

\section{Irreducible Unitary Representations of $\tilde P$}
Let $G$ be the symplectic group $Sp_{2n}(\mb R)$ with {\bf $n \geq 2$} and $P$ be the maximal parabolic subgroup preserving a one dimensional isotropic subspace $\mathbb R e_1$. Let $\tilde G$ be the universal covering of $G$.  For simplicity, let $\mathbb Z$ be the preimage of the identity. The group $P$ has a Langlands decomposition $GL_1(\mathbb R) Sp_{2n-2}(\mb R) H_{n-1}$ where $H_{n-1}$ is the Heisenberg group. $P$ is a semidirect product of
$GL_1(\mb R) \times  Sp_{2n-2}(\mb R)$ and $H_{n-1}$. $GL_1(\mb R)$ can be further decomposed as $\mathbb Z_2 A $ with $A \cong \mb R^+$.
\begin{lem}  $\tilde A \cong \mathbb Z \tilde A_0$ where $\tilde A_0$ is the identity component of $\tilde A$ which can be identified with $A$. So we will write $\tilde A= \mathbb Z A$. In addition,
$\widetilde{GL}_1(\mathbb R) \cong (\frac{1}{2} \mathbb Z ) A$. Lastly
$$ \tilde P \cong  (\frac{1}{2} \mathbb Z)  A \widetilde{Sp}_{2n-2}(\mb R)  H_{n-1}/ \mathbb Z.$$ 
So $\pi_0(\tilde P)=\mathbb Z_2$. Here $\widetilde{GL}_1(\mb R) \cap \widetilde{Sp}_{2n-2}(\mb R)= \mathbb Z $.
\end{lem}
Notice that the adjoint action of $\widetilde{GL}_1(\mathbb R)$ on $H_{n-1}$ descends into the adjoint action of $GL_1(\mb R)$ on $H_{n-1}$ and the adjoint action of $\widetilde{Sp}_{2n-2}(\mb R)$ descends into the adjoint action of $Sp_{2n-2}(\mb R)$ on $H_{n-1}$.\\
\\
Suppose that $\lambda$ is real and $\lambda \neq 0$. Let $\rho_{\lambda}$ be the unique irreducible unitary representation of $H_{n-1}$ with central character $\exp i \lambda t$. The adjoint action of $GL_1(\mb R)$ on $H_{n-1}$ induces an action of $GL_1(\mb R)$ on $\hat{H}_{n-1}$. In particular, $\pm 1 \in GL_1(\mathbb R)$ perserve $\rho_{\lambda}$ and
$$a \in GL_1(\mb R): \rho_{\lambda} \rightarrow \rho_{a^2 \lambda}.$$
By Mackey analysis, there are two irreducible unitary representations of $A H_{n-1}$:
$$\rho_{+}=\int_{\lambda \in \mathbb R^+} \rho_{\lambda} d \lambda, \qquad \rho_{-}=\int_{\lambda \in \mathbb R^{-}}
\rho_{\lambda} d \lambda.$$
These are the only irreducible unitary representations with $\rho|_{C(H_{n-1})} \neq I$, the identity.
Now $\pm 1 \in GL_1(\mb R)$ preserves each $\rho_{\lambda}$. By Stone-Von Neumannn Theorem, $\pm 1$ acts on each $\rho_{\lambda}$ projectively. In this situation, it is easy to make $\pm 1$ act on $\rho_{\lambda}$ directly. There is no obstruction to lift the projective action of $\pm 1$ on $\rho_{\lambda}$. Using the Schr\"odinger model,
$-1$ acts on the odd functions by $-1$ and on the even functions by $+1$. Let us include the actions of $\mathbb Z_2 \subseteq GL_1(\mb R)$ in the model $\rho_{\lambda}$, consequently in $\rho_{\pm}$. Now again, by Mackey analysis,
there are four irreducible unitary representations of $GL_1(\mb R) H_{n-1}$ on which $C(H_{n-1})$ acts nontrivially, namely
$$\rho_{\pm} \otimes {\rm sgn}, \qquad \rho_{\pm}.$$
The difference between the former and the latter is a little subtle. One way to tell the difference is that
$(\rho_{\pm} \otimes {\rm sgn})(-1)$ acts on the even functions by $-1$ while $\rho_{\pm}(-1)$ acts on the even functions by identity. \\
\\
Now consider $\widetilde{GL}_1(\mb R) H_{n-1}$. The representation $\rho_{\pm}$ can be regarded as a representation of $\widetilde{GL}_1(\mb R) H_{n-1}$. 
\begin{lem} Identify $\widetilde{GL}_1(\mb R) H_{n-1}$ with $(\frac{1}{2} \mathbb Z ) A  H_{n-1}$. Then the irreducible unitary representations on which $C(H_{n-1})$ act nontrivially are all of the form
$$\rho_{\pm} \otimes \chi_t \mid t \in [0, 1) $$
with $\chi_t(m) =\exp 4 \pi i m t$ for $ m \in \frac{1}{2}\mathbb Z$.
\end{lem}
For $GL_1(\mb R) H_{n-1}$, $t=0, \frac{1}{2}$ because $\chi_t(\mathbb Z)=1$.\\
\\
Now let us consider $\widetilde{Sp}_{2n-2}(\mb R)$. This group preserves $\rho_{\lambda}$. Again, by the Stone-Von Neumann Theorem, $\widetilde{Sp}_{2n-2}(\mb R)$ acts on $\rho_{\lambda}$ projectively. Since $\widetilde{Sp}_{2n-2}(\mb R)$ is already simply connected, one obtains a group action of $\widetilde{Sp}_{2n-2}(\mb R)$ on $\rho_{\lambda}$. By a theorem of Segal-Shale-Weil, $\widetilde{Sp}_{2n-2}(\mb R)$ action on $\rho_{\lambda}$
descends into an action of $Mp_{2n-2}(\mb R)$. Simply put, $m \in \mathbb Z$ acts by $(-1)^m=\exp i m \pi$ on $\rho_{\lambda}$. We can now extend $\rho_{\pm}$ to include the action of $\widetilde{Sp}_{2n-2}(\mb R)$. \\
\\
 {\bf From now on
$\rho_{\pm}$ will be  representations of $\tilde P$}.\\
\begin{thm}[\cite{howe}] \label{cla} Irreducible unitary representations of $\tilde P$ on which $C(H_{n-1})$ acts nontrivially are of the form
$$[\rho_{\pm} \otimes \chi_t] \otimes \tau$$
with $\tau$ an irreducible unitary representation of $\widetilde{Sp}_{2n-2}(\mb R)$ such that
$\tau(m)= \chi_t(m)$ for any $m \in \mathbb Z$. In addition, two such representations are equivalent if and only
if all the parameters $(\pm, t, \tau)$ are the same.
\end{thm}
I shall make some remarks here. First, $\tau$ is extended to a representation of $\tilde P$, trivially on $ A H_{n-1}$, and trivially on the component group. Second, $\chi_t \otimes \tau$ is a twisted tensor product in the sense that the action of $\mathbb Z$ commutes with the tensor.
So $\chi_t(m) \otimes \tau(n)=\chi_t(n) \otimes \tau(m)$ for any $m,n \in \mathbb Z$. For group $P$, $\rho_{\pm}(m) \otimes \chi_{t}(m) \otimes \tau(0)$ must be the identity for every $m \in \mathbb Z \subseteq \tilde G$. So $(-1)^m \exp 4 \pi m t =(-1)^m  \tau(m) = 1$. For an irreducible unitary representation of $P$ on which $C(H_{n-1})$ acts nontrivially,  $t = \frac{1}{4}, \frac{3}{4}$ and $\tau$ is a genuine unitary representation of $Mp_{2n-2}(\mb R)$. \\
\\
The proof is straight forward by applying the Mackey analysis. Observe that the subgroup of $\tilde P$ that preserves $\rho_{\pm}$ is 
$\frac{1}{2} \mathbb Z \times_{\mathbb Z} \widetilde{Sp}_{2n-2}(\mb R)$.  $\chi_t \otimes \tau$ parametrizes
the equivalence classes of irreducible unitary representations of this subgroup.

\section{Irreducible Unitary Representations of $\widetilde{Sp}_2(\mb R)$}
Throughout this section $G=Sp_{2}(\mb R)$ and $P=MAN$ where $M \cong \mathbb Z_2$ and $A \cong \mathbb R^+$ and $N \cong \mathbb R$. What we have proved in the last section needs to be modified. Since the unitary dual of $\tilde G$ is known (\cite{puk}, \cite{ht}), we will analyze $\pi|_{\tilde P}$ in detail. The results in this section must have been known to the experts. They will be used in the next section to analyze higher rank case. \\
\\
Fix the standard maximal compact group $SO(2)$. We parametrize it by the angle of rotation {\it counterclockwise}.
There are essentally four classes of irreducible unitary representations of $\tilde G$ (see ~\cite{puk} \cite{knapp}):
\begin{enumerate}
\item the trivial representation $\mathds 1$;
\item unitary principal series $I(\epsilon, s)$ where $\epsilon \in [0, 1)$ and $s \in i \mathbb R$ ( we exclude
$\epsilon=\frac{1}{2}, s=0$);
\item complementary series $C(\epsilon, s)$ ($\epsilon \in [0, 1)$, $ s \in (0, | 1-  2 \epsilon|)$);
\item Highest weight modules $D_{l }^{-} (l > 0)$ and lowest weight modules $D_{l}^{+}(l >0)$.
\end{enumerate}
Let $P$ be the standard upper triangular parabolic subgroup of $G$. Let $N$ be the nilradical of $P$. Then the identity component $P_0$ has two irreducible unitary representations on which $N$ acts nontrivially, namely $\rho_{+}$ and $\rho_{-}$. $\rho_{+}|_{N}$ is supported on $\mathbb R^+ \subseteq \hat N$, and $\rho_{-}|_{N}$ is supported on $\mathbb R^{-} \subseteq \hat N$. Now the center of $\tilde G$ can be identified with $\frac{1}{2} \mathbb Z$.  Identifying $\tilde P$ with $\frac{1}{2} \mathbb Z A N$, $\rho_{\pm}$ can be extended to a representation of $\tilde P$ by identify $P_0$ with $\tilde P/ \frac{1}{2} \mathbb Z$. Then elements of $\hat{\tilde P}$ are parametrized by $(\pm, t)$ with $t \in [0, 1)$. More precisely, every irreducible unitary representation of $\tilde P$ is equivalent to $\rho_{\pm} \otimes \chi_t$. Here $\chi_t(man)=\exp 4 \pi i m t$ with $m \in \frac{1}{2} \mathbb Z$, $a \in A$ and $n \in N$. Notice that $\chi_t$ also defines a central character of $\tilde G$. In our setting, the representations with even weights have central character $\chi_0$; the representations with odd weights have central character $\chi_{\frac{1}{2}}$. Let $\chi_1=\chi_0$.\\
\\
The following theorem gives the structure of the restriction of irreducible unitary representations of $\tilde G$ to $\tilde P$. 

\begin{thm}\label{real} Let $\lfloor \frac{l}{2} \rfloor$ be the largest integer less or equal to $\frac{l}{2}$.
\begin{enumerate}
\item $I(\epsilon, s )|_{\tilde P} \cong (\rho_{+} \oplus \rho_{-}) \otimes \chi_{\epsilon};$
\item $C(\epsilon, s)|_{\tilde P} \cong (\rho_{+} \oplus \rho_{-}) \otimes \chi_{\epsilon};$
\item $D_{l}^+ |_{\tilde P} \cong \rho_{+} \otimes \chi_{\frac{l}{2}- \lfloor \frac{l}{2} \rfloor};$
\item $D_{l}^{-}|_{\tilde P} \cong \rho_{-} \otimes \chi_{1- \frac{l}{2} +\lfloor \frac{l}{2} \rfloor}.$

\end{enumerate}
\end{thm}

The results here are obviously known to the experts. We will provide an elementary proof.\\

Proof: The central character of each $\pi \in \hat{\tilde G}$ can be computed easily.  Using the noncompact model, $I(\epsilon, s)$ can be modeled on $L^2(N)$ with $N$ act as translations. Hence (1) is proved.  To prove $(2), (3), (4)$, it suffices to show that the Fourier transform of the matrix coefficients of $\pi$ restricted to $N$ has the desired support.\\
\\
To show (2), let $(\, , \, )$ be the inner product of $C(\epsilon, s)$ and $(\, , \,)_{\rm Inv}$ be the natural complex linear pairing between the induced representations $I^{\infty}(\epsilon, -s)$ and $I^{\infty}(\epsilon, s)$. For smooth vectors $\phi, \psi \in I^{\infty}(\epsilon, s)$, we have
$$(\phi, \psi)=(A(\epsilon, s) \phi, \psi)_{\rm Inv}$$
where $A(\epsilon, s)$ is the intertwining operator defined over smooth vectors. In addition, $A(\epsilon, s)$ defines a bijection between $I^{\infty}(\epsilon, s)$ and $I^{\infty}(\epsilon, -s)$.
Using the noncompact model,  for every $n \in N$ as an additive group, we have
$$(C(\epsilon, s)(n) \phi, \psi)=(A(\epsilon, s) \phi, I(\epsilon, s)(-n) \psi)_{\rm Inv}=\int_{N} A(\epsilon, s)\phi(x) \overline{\psi}(x+n) d x .$$
Now $C^{\infty}_c(N) \subseteq I^{\infty}(\epsilon, \pm s)$. We choose $\phi$ and $\psi$ so that 
$$A(\epsilon, s) \phi \in C^{\infty}_c(N), \qquad \psi \in C^{\infty}_c(N) .$$
So $(C(\epsilon, s)(n) \phi, \psi)$ becomes the convolution of two smooth and compactly supported functions. Its Fourier transform can be made to be supported on $\hat N$, upon proper choices of $\phi$ and $\psi$.
So again we have $C(\epsilon, s)|_{\tilde P} \cong (\rho_{+} \oplus \rho_{-}) \otimes \chi_{\epsilon}$.\\
\\
We will now prove (3). (4) follows immediately from (3). Let $l >0$. Notice that $D_{l}^{+} $ is a subquotient of $I(\epsilon, l-1)$ and of $I(\epsilon, 1-l)$ with $\epsilon=\frac{l}{2}-\lfloor \frac{l}{2} \rfloor$. Let $v_{l+p}$ be of weight $l+p$ in $D_l^{+}$ where $p$ is a nonnegative even integer.   We stick with the noncompact picture. Let 
$$\phi_{l+p}(x)=(\frac{1}{\sqrt{x^2+1}})^l (\frac{1-x i}{1+xi})^{\frac{l+p}{2}}=\frac{(1+x^2)^{\frac{p}{2}}}{(1+ xi)^{l+p}}$$ be a function in the noncompact model of $I(\epsilon, l-1)$. Here we choose the standard $\arg$ function between $-\frac{\pi}{2}$ and $\frac{\pi}{2}$ to define $(1 \pm x i)^{l}$ if $l$ is not an integer. 
Let 
$$\psi_{l+q}=(\frac{1}{\sqrt{x^2+1}})^{2-l}(\frac{1-x i}{1+xi})^{\frac{l+q}{2}}=\frac{(1- xi)^{l-1+\frac{q}{2}}}{(1+xi)^{1+\frac{q}{2}}}$$ be a function in the noncompact model of $I(\epsilon, 1-l)$. Then 
$$(D_{l}^{+}(g) v_{l+p}, v_{l+q})= C (I(\epsilon, l-1)(g) \phi_{l+p}, \psi_{l+q})_{\rm Inv} \qquad (g \in \tilde G, p \geq 0, q \leq 0).$$
In particular, for $n \in N \cong \mathbb R$, $(D_{l}^{+}(n) v_{l+p}, v_{l+q})= C \int_{N} \phi_{l+p}(x-n) \psi_{l+q}(x) d x$. \\
\\
Recall that 
$$ \frac{1}{(1+x i)^{l+p}}= C_1 \int_{\mathbb R^+} (\exp - \xi) \, \xi^{l+p-1} (\exp  (-i x \xi )) \, d \xi $$
(See Ch. 8.3 \cite{fo}). It follows  by Fourier analysis that
$$\phi_{l+p}(x)= \frac{(1+x^2)^{\frac{p}{2}}}{(1+ x i)^{l+p}}= \int_{\mathbb R^+} W(\xi) \exp  (-i x \xi ) \, d \xi $$
where $W(\xi)$ is a linear combination of derivatives of $(\exp - \xi) \, \xi^{l+\frac{p}{2}-1}|_{\mathbb R^+}$ and $W(\xi)$ is supported on $\mathbb R^+$. Therefore
$$\phi_{l+p}(x-n)=\int_{\mathbb R^+} W(\xi) \exp  (-i x \xi ) \exp (i n \xi) \, d \xi .$$
 Even though  $\psi_{l+q}(x)$ are not in $L^1$,  $\frac{d^{\lceil  l \rceil+3}}{d x^{\lceil  l \rceil+3}} \psi_{l+q}(x)$ will be in $L^1$. So Fourier transform of $\psi_{l+q}(x)$ will be a $C_0$ function on $\mathbb R$ multiplied by a monomial of $\xi$. It is easy to see that the Fourier transform of $(D_{l}^{+}(n) v_{l+p}, v_{l+q})$ is supported on $\mathbb R^+$. Since $\{ v_{l+p} \mid p>0 \}$ is an orthogonal basis for $D_{l}^+$, 
$D_{l}^+|_{N}$ is supported on $\mathbb R^+$.  $\Box$\\
\\
We shall remark that Theorem ~\ref{real} depends on the parametrization of  $K$ and $N$. If one chooses the opposite parabolic $\overline P$, then the statements in $(3)$ and $(4)$ will change. There may be other proofs to Theorem ~\ref{real}.  The proof of (3) and (4) we give here is more self-contained.

\section{Restriction of Unitary Representations and Irreducibility}
Now let $\pi$ be a nontrivial  irreducible unitary representation of $\tilde G$. Then $\pi|_{\tilde P}$ can be decomposed into a direct integral of $\rho_{\pm} \otimes \chi_t \otimes \tau$. In particular, one can write
$$\pi  \cong \rho_{+} \otimes \tau_{+}(\pi) \oplus \rho_{-} \otimes \tau_{-}(\pi) \oplus \tau_0.$$
Here $\tau_{+}(\pi)$ and $\tau_{-}(\pi)$ are unitary representations of $\widetilde{GL}_1(\mb R) \widetilde{Sp}_{2n-2}(\mb R)$. The following theorem says that $\tau_0$ does not occur.\\

\begin{thm}[See \cite{howe}]\label{rho} Let $\pi$ be a nontrivial irreducible unitary representation of $\tilde G$. Then there exist two unitary representations $\tau_{+}(\pi)$ and $\tau_{-}(\pi)$ of $\widetilde{GL}_1(\mb R) \widetilde{Sp}_{2n-2}(\mb R)$ such that
$$\pi|_{\tilde P}  \cong \rho_{+} \otimes \tau_{+}(\pi) \oplus \rho_{-} \otimes \tau_{-}(\pi).$$
\end{thm}
Notice that one of $\tau_{\pm}(\pi)$ could be zero. This theorem was proved by Howe in \cite{howe} Pg. 249 for the metaplectic group.\\
\\
Proof: We will have to prove that $\pi|_{\tilde P}$ does not have any subrepresentation on which $C(H_{n-1})$ acts trivially.  Suppose otherwise. let  $v$ be a nonzero vector fixed by $C(H_{n-1})$. Let $\tilde G_0$ be the subgroup of $\tilde G$ that commutes with $\widetilde{Sp}_{2n-2}(\mb R)$.
So $\tilde G_0 \cong \widetilde{Sp}_2(\mb R)$. Notice that $C(H_{n-1}) \subseteq \tilde G_0$ and $A \subseteq \tilde G_0$.  Let $\mathcal H$ be the Hilbert space spanned by $\pi(\tilde g_0) v$ for $\tilde g_0 \in \tilde G_0$.
Clearly, $\mc H$ decomposes into a direct integral of irreducible unitary representations of $\tilde G_0$ on which
$\mathbb Z$ acts as a character. Indeed, all factorial representations of $\tilde G_0$ are direct sum of irreducible
representations. Now let $$v=\int_{\widehat{\tilde G}_0} v_s d \mu(s)$$
 where $v_s \in \mc H_s \otimes V_s$, $\mc H_s \in \hat{\tilde G}_0$ and $\dim V_s=m(\mc H_{s}, \mc H)$. Then $C(H_{n-1})$ must fix $v_s$ for almost all $s$ with respect to $\mu$. If $\mc H_s$ is not trivial, $\mc H_s \otimes V_s$ has no vector fixed by $C(H_{n-1})$. Hence, $\mc H$ must be a direct sum of the trivial representation of $\tilde G_0$. In particular, $\pi$ must descend to a representation of $G$. The matrix coefficent $g \rightarrow (\pi(g)v, v)$ violates the Howe-Moore vanishing Theorem (~\cite{hm}). We reach a contradiction. $\Box$ \\
 \\
Now we can prove one of our main results.
 
\begin{thm}\label{main2} Let $\pi$ be a nontrivial irreducible unitary representation of $\tilde G$ such that $\mathbb Z$ acts by $\exp 2 \pi m t$ ($ \forall \, m \in \mathbb Z$) for a fixed $t \in [0, 1)$. Suppose that $\pi|_{\tilde P}$ is irreducible. Then $\pi$ must be a highest weight module or a lowest weight module.
\end{thm}
Proof: Let us fix the standard maximal compact group $U(n) \subseteq G$. Then $U(n) \cap Sp_{2n-2}(\mb R)=U(n-1)$. As usual, the complexified Lie algebra $\mathfrak g_{\mathbb C}$ decomposes into a direct sum
$$ \mathfrak k_{\mathbb C} \oplus \f p^+ \oplus \f p^-.$$
Suppose that $\pi|_{\tilde P}$ is irreducible. By the last Theorem, either $\tau_{+}(\pi)$ or $\tau_{-}(\pi)$ must be zero. Without loss of generality, suppose that $\pi|_{\tilde P} \cong \rho_{+} \otimes \tau_{+}(\pi)$.  Notice that $\rho_{+}|_{C(H_{n-1})}$ is supported on $\mathbb R^+$. So $\pi|_{C(H_{n-1})}$ must also be supported on $\mathbb R^+$.  Let $\tilde G_0$ be the subgroup of $\tilde G$ that commutes with $\widetilde{Sp}_{2n-2}(\mb R)$. Consider the restriction $\pi|_{\tilde G_0}$. $\pi|_{\tilde G_0}$ can be decomposed into a direct integral of irreducible unitary representations with multiplicities. By Theorem \ref{real}, among the irreducible unitary representations of $\tilde G_0$, only the lowest weight modules are supported on $\mathbb R^+ \subset \widehat{C(H_{n-1})}$. Hence only the lowest weight modules occur in the direct integral decomposition of $\pi|_{\tilde G_0}$. \\
\\
Now $\pi|_{\tilde G_0}$ is a direct integral of lowest weight modules. Let $U(1)=G_0 \cap U(n)$. Then $\pi|_{\tilde U(1)}$ can only have positive weights. Fix a maximal torus $ T \supseteq U(1)$ in $U(n)$.  $\pi|_{\tilde{T}}$ can only have positive weights, since the weight space 
of $\pi|_{\tilde{T}}$ is invariant under the Weyl group of $U(n)$.  Let $v_{\lambda}$ be a vector with weight $\lambda$ such that $\sum \lambda_i$ is minimal among all possible weights occuring in $\pi|_{\tilde T}$. Notice that the set of all possible $\sum \lambda_i $ is a discrete set in $\mathbb R^+$. So a minimal $\sum \lambda_i$ must exist. Now $\pi|_{\tilde U(n)}$ must contain an irreducible representation $V_{\mu}$ with $\sum \mu_i=\sum \lambda_i$. Clearly, $\f p^-$ act on $V_{\mu}$ by zero. So the module generated by $V_{\mu}$ must have a lowest weight module as its quotient. Since $\pi$ is already unitary and irreducible, $\pi$ must be a lowest weight module. Now we have shown that $\pi$ is a unitary lowest weight module. $\Box$ \\
\\
We shall remark that the last paragraph is true even one assumes that the weights for $\pi|_{\tilde U(1)}$ is bounded from below. 

\section{Some Criterions for Lowest Weight Modules}

In this section, we give some characterization of lowest weight modules in terms of their restrictions on certain subgroups. Some of them are well-known to the experts. Let us fix a complex structure and an inner product $( \, , \,)$ on the symplectic space $\mathbb R^{2n}$ such that the symplectic form coincides with the imaginary part of $(\, , \,)$. Let $e_1, e_2, \ldots, e_n$ be the standard basis over $\mathbb C$. Let $G=Sp_{2n}(\mb R)$. Let $P$ be the subgroup preserving $\mathbb R e_1$. Then $P=GL_1(\mb R) Sp_{2n-2}(\mb R) H_{n-1}$. Let $U(n)$ be the subgroup preserving $(\, , \,)$.\\

\begin{thm}~\label{equ}
Let $T$ be the standard maximal torus in the maximal compact group $U(n)$ of $G$. Let $Sp_2(\mb R)$ be the subgroup of $G$ preserving $\mathbb C e_1$ and acting on the complex linear span of $\{ e_2, e_3, \ldots, e_{n} \}$ by identity. Let $U(1)=Sp_2(\mb R) \cap U(n)$. Let $\pi$ be a nontrivial irreducible unitary representation of $\tilde G$. Let $Q$ be a maximal parabolic subgroup of $G$. Let $N$ be its nilradical and $ZN$ be the center of $N$. Suppose that $ZN \supseteq C(H_{n-1})$.
The following are equivalent:
\begin{enumerate}
\item  
$\pi|_{ZN}$ is supported on a subset of the positive semidefinite cone of $\widehat{ZN}$, regarded as the space of symmetric matrices;
\item $\pi|_{C(H_{n-1})}$ is supported on  $\mathbb R^+ \subseteq \widehat{C(H_{n-1})}$;
\item $\pi|_{\widetilde{Sp}_{2}(\mb R)}$ decomposes into a direct integral of lowest weight modules;
\item $\pi|_{\tilde U(1)}$ only has positive weights;
\item $\pi|_{\tilde T}$ only has positive weights;
\item the weights of $\pi|_{\tilde U(1)}$ is bounded from below;
\item there is an integer $k$ such that every weight $\lambda$ of $\pi|_{\tilde T}$ satisfies $\lambda_i \geq k$ for every $i$;
\item $\pi$ is a unitary lowest weight module of $\tilde G$;
\item If $n=1$, $\pi$ is a unitary lowest weight module; if $n \geq 2$, $\pi|_{\tilde P}$ decomposes into $\rho_{+} \otimes \tau_{+}(\pi)$ and the weights  of $\tau_{+}(\pi)|_{\tilde T(1)}$ are bounded from below. Here $T(1)$ is the one dimensional compact torus in $U(n) \cap Sp_{2n-2}(\mb R)$ fixing all vectors in the complex span of $\{ e_1, e_3, \ldots, e_n \}$.
\end{enumerate}
\end{thm}
Proof:  When $n=1$, our theorem follows from Theorem \ref{real}. Suppose now $n \geq 2$. \\
\\
By ~\cite{howe}, $\pi|_{ZN}$ is supported on $GL$-orbits on $\widehat{ZN}$. $(1) \leftrightarrow (2)$ is a matter of matrix analysis. $(2) \rightarrow (3) \rightarrow (4) $ is proved in Theorem ~\ref{main2}.  $(4) \rightarrow (3) \rightarrow (2)$ is easier than the other direction. So $(1), (2), (3)$ and $(4)$ are equivalent. \\
\\
$(4) \rightarrow (5) \rightarrow (6) \rightarrow (7)$ is trivial. $(7) \rightarrow (8)$ follows as in Theorem ~\ref{main2}. $(8) \rightarrow (6)$ is also obvious. So $(8), (7), (6)$ are equivalent. \\
\\
To prove $(8) \rightarrow (9)$, suppose $\pi$ is a nontrivial unitary lowest weight module. By Cor. ~\ref{rho}, $\pi|_{\tilde P} \cong \rho_{+} \otimes \tau_{+}(\pi) \oplus \rho_{-} \otimes \tau_{-}(\pi)$. If $\tau_{-}(\pi) \neq 0$, fix a $\tilde T(1)$-eigenvector $v$ with weight $\lambda$. By tensoring with vectors in $\rho_{-}$, we obtain $\tilde T(1)$-eigenvector with arbitarily low weight. By Weyl group action, we obtain $\tilde U(1)$-eigenvectors with arbitrarily low weight. This contradicts $(6)$. So
$\pi \cong \rho_{+} \otimes \tau_{+}(\pi)$. Similarly, weights for $\tau_{+}(\pi)|_{\tilde T(1)}$ must also be bounded below. Hence $(8) \rightarrow (9)$.\\
\\ 
Suppose that $(9)$ holds. So  the weights  of $\tau_{+}(\pi)|_{\tilde T(1)}$ are bounded from below. Let $G_{(2)}$ be the subgroup of $Sp_{2n-2}(\mb R)$ that fixes all vectors in the complex linear span of $\{e_1, e_3, e_4, \ldots e_n \}$. So $T(1)$ is maximal compact in $G_{(2)}$ and
$\tau_{+}(\pi)|_{\widetilde{G_{(2)}}}$ must be a direct integral of lowest weight modules with multiplicites. So the weights of $\tau_{+}(\pi)|_{\tilde T(1)}$ must all be positive. Hence the weights of $\pi|_{\tilde T(1)}$ must all be positive. Due to the action of the Weyl group of $U(n)$, $\pi|_{\tilde U(1)}$ only has positive weights. So $(9) \rightarrow (4)$. We have proved $(9) \rightarrow (4) \rightarrow (5) \rightarrow (6) \rightarrow (7) \rightarrow (8) \rightarrow (9)$.  $\Box$ \\
\\
 We shall remark here that
the parametrization of unitary highest or lowest weight modules is already known, due to the work of Enright-Howe-Wallach (~\cite{ehw}). It is perhaps easy to go further to derive more properties of unitary lowest or highest weight modules from Theorem ~\ref{equ}.
\begin{cor} Let $\pi$ be an irreducible unitary lowest weight module of $\tilde G$. Then $\pi|_{\tilde P} \cong \rho_{+} \otimes \tau_{+}(\pi)$ and $\tau_{+}(\pi)$ decomposes into a direct sum of lowest weight modules of $\widetilde{Sp}_{2n-2}(\mb R)$.
\end{cor}
Proof: By Theorem \ref{equ} (9), the weights of $\tau_{+}(\pi)|_{\tilde U(1)}$ are bounded from below. By Theorem \ref{equ} (4), $\tau_{+}(\pi)|_{\tilde U(1)}$ only has positive weights. By Theorem \ref{equ} (7) (8), there is a lowest weight subrepresentation
of $\widetilde{Sp}_{2n-2}(\mb R)$ in $\tau_{+}(\pi)$. Consider the orthogonal complement. If it is nonzero, then there is another lowest weight subrepresentation. This process can continue and it will end in countable time due to the fact that $\pi$ has a countable basis. We now obtain a discrete decomposition. $\Box$

\section{The group $U(p,q)$}
Suppose $p \geq q \geq 1$ and $p+q \geq 3$. Let $U(p,q)$ be the group that preserve a Hermitian form $(\, , \,)$ on $\mathbb C^{p+q}$ with signature $(p,q)$. Let $P$ be a maximal parabolic subgroup that preserves a one dimensional isotropic subspace. Then $P$ can be identified with $GL_1 (\mb C) U(p-1, q-1) H_{p+q-2}$. Here $H_{p+q-2}$ are parametrized by $(t \in \mathbb R, u \in \mathbb C^{p+q-2})$ . The adjoint action of $g \in U(p-1, q-1)$ on $H_{p+q-2}$ leaves $t$ fixed and operates on $u$ as the left multiplication. The adjoint action of $a \in GL_1(\mb C)$ on $H_{p+q-2}$ dilates $t$ to $\|a\|^2 t$ and operates on $u$ as scalar multiplication. 
Write $GL_1(\mb C)= A U(1)$ where $A =\mathbb R^+$. \\
\\
Let $\tilde G=\{ (g, t) \mid g \in U(p,q), \exp 2 \pi i t =\det g \}$ be an infinite covering of $G$. For any subgroup $H$ of $G$, let $\tilde H$ be the preimage of $H$.
Let $\rho_{\pm}$ be the irreducible unitary representations of $A H_{p+q-2}$ we defined earlier. These are the only irreducible unitary representations on which $C(H_{p+q-2})$ acts nontrivially. Extend $\rho_{\pm}$  to a unitary representation of $\tilde P$. Notice that $\rho_{+}|_{\tilde{U}(1) \tilde{U}(p-1, q-1)}$ decomposes as follows
$$ [ \oplus_{n \in \mathbb Z} {\det \,}^{n+\frac{p-q}{2}} \otimes \theta({\det \,}^{n+\frac{p-q}{2}}) ] \otimes \mathbb C^{\infty}.$$
Here ${\det \,}^{n+\frac{p-q}{2}}$ is a character of $\tilde U(1)$, $\theta({\det \,}^{n+\frac{p-q}{2}})$ is the theta lift of ${\det \,}^{n+\frac{p-q}{2}}$ with respect to $(U(1), U(p-1, q-1))$ (see \cite{kv} \cite{h} \cite{paul}) and $\mathbb C^{\infty}$ records the multiplicity.\\

\begin{thm}~\label{main1u}
Let $\pi$ be an infinite dimensional irreducible unitary representation of $\tilde G$. Then there are two unitary representation $\tau_{\pm}(\pi)$ of $\tilde{U}(p-1, q-1)$ such that
$$\pi|_{\tilde P} \cong \rho_{+} \otimes \tau_{+}(\pi) \oplus \rho_{-} \otimes \tau_{-}(\pi).$$
Let $U(1)$ be diagonally embedded in $U(1,1) \subseteq U(p,q)$. Then
\begin{equation}
\begin{split}
\pi|_{\tilde{U}(1) \tilde{U}(p-1, q-1)} \cong &  [ \oplus_{n \in \mathbb Z} {\det \,}^{n+\frac{p-q}{2}} \otimes \theta({\det \,}^{n+\frac{p-q}{2}})] \otimes \tau_{+}(\pi) \otimes \mathbb C^{\infty} \\
& \oplus  [ \oplus_{n \in \mathbb Z} {\det \,}^{n+\frac{p-q}{2}} \otimes \theta^{\prime}({\det \,}^{n+\frac{p-q}{2}})] \otimes \tau_{-}(\pi) \otimes \mathbb C^{\infty}.
\end{split}
\end{equation}
Here $\theta^{\prime}$ refers to the theta lifts with respect to the contragredient oscillator representation.
\end{thm}
Notice that one of the $\tau_{\pm}(\pi)$ could be zero. 
Similarly, we can prove \\

\begin{thm}~\label{main2u} Let $\pi$ be an infinite dimensional irreducible unitary representation of $\tilde G$. Suppose that $\pi|_{\tilde P}$ is irreducible. Then $\pi$ must be a highest weight module or a lowest weight module.
\end{thm}
Essentially, $SU(1,1)$ in $SU(p,q)$ will play the role of $Sp_2(\mb R)$ in $Sp_{2n}(\mb R)$. The proof is omitted here.
\section{The group $O^*(2n)$}
Let $n \geq 3$. Let $O^*(2n)$ be the group of isometry preserving a nondegenerate skew-Hermitian form on $\mathbb H^n$. Let $P$ be the maximal parabolic subgroup preserving a 1 dimensional isotropic subspace. Then $P$ can be identified with
$GL_1(\mathbb H)O^*(2n-4)H_{2n-4}$ where $H_{2n-4}$ is a Heisenberg group parametrized by 
$$(t \in \mathbb R, u \in \mathbb H^{n-2}).$$
$GL_1(\mathbb H)$ can be further decomposed as $Sp(1) A$ where $A$ is the center of $GL_1(\mathbb H)$.
The adjoint action of $ a \in A$ on $H_{2n-4}$ is given by
$$ (t, u) \rightarrow (a^2 t, a u).$$
The adjoint action of $O^*(2n-4)$ on $H_{2n-4}$ is the left multiplication on $u$. The adjoint action of $k \in Sp(1)$ on $H_{2n-4}$ is the right multiplication. Clearly, $Sp(1) \times O^*(2n-4)$ action preserves 
the real part of canonical skew-Hermitian form on $\mathbb H^{n-2}$. $(Sp(1), O^*(2n-4))$ becomes a dual reductive pair (See \cite{h}).\\
\\
Now let $\rho_{\pm}$ be the two irreducible unitary representations of $AH_{2n-4}$ on which $C(H_{2n-4})$ acts nontrivially. $\rho_{\pm}$ extends to irreducible unitary representations of the linear group $P$. In particular, $\rho_{+}|_{Sp(1) O^*(2n-4)}$ decomposes according to the theta correspondence with infinite multiplicity:
$$ \oplus_{\sigma \in \widehat{Sp(1)}} \sigma \otimes \theta(\sigma) \otimes \mathbb C^{\infty}.$$
Similarly, 
$$\rho_{-}|_{Sp(1) O^*(2n-4)} \cong \oplus_{\sigma \in \widehat{Sp(1)}} \sigma \otimes \theta^{\prime}(\sigma) \otimes \mathbb C^{\infty}.$$
Here $\theta^{\prime}$ is the theta correspondence with respect to the contragredient oscillator representation.
$\widehat{Sp(1)}$ is parametrized by $\mathbb N$.\\

\begin{thm}~\label{main1o}
Let $\pi$ be a nontrivial irreducible unitary representation of $O^*(2n)$. Then there exists two unitary representations $\tau_{\pm}(\pi)$ of $O^*(2n-4)$ such that 
$$\pi|_{P} \cong \rho_{+} \otimes \tau_{+}(\pi) \oplus \rho_{-} \otimes \tau_{-}(\pi).$$
In addition,
$$\pi|_{Sp(1) O^*(2n-4)} \cong \{ \oplus_{\sigma \in \widehat{Sp(1)}} [ \sigma \otimes \theta(\sigma)] \otimes \tau_{+}(\pi) \otimes \mathbb C^{\infty} \} \oplus \{ \oplus_{\sigma \in \widehat{Sp(1)}} [ \sigma \otimes \theta^{\prime} (\sigma)] \otimes \tau_{-}(\pi) \otimes \mathbb C^{\infty} \}$$
\end{thm}
One of $\tau_{\pm}(\pi)$ could be zero. The theorem for the universal covering of $O^*(2n-4)$ is left to the reader. Similarly, we have \\

\begin{thm}~\label{main2o} Let $\pi$ be a nontrivial irreducible unitary representation of $O^*(2n)$. If $\pi|_P$ is irreducible, then $\pi$ must be a highest weight module or lowest weight module.
\end{thm}
Notice that the group $O^*(4)$ contains a noncompact factor $SU(1,1) \cong SL_2(\mb R)$. The proof is essentially the same as in Theorem \ref{main2}.

\end{document}